\documentclass[12pt]{article}
\usepackage{amsfonts,amsthm}
\usepackage{epsfig}
\usepackage{makeidx}
\usepackage{graphicx}

\makeindex
\newenvironment{namelist}[1]{%
\begin{list}{}
{

\settowidth{\labelwidth}{#1}
\setlength{\leftmargin}{1.1\labelwidth}
}
}{%
\end{list}}
\newcommand{\ncom}{\newcommand}
\ncom{\ul}{\underline}
\ncom{\beq}{\begin{equation}}
\ncom{\eeq}{\end{equation}}
\ncom{\bea}{\begin{eqnarray*}}
\ncom{\eea}{\end{eqnarray*}}
\ncom{\beqa}{\begin{eqnarray}}
\ncom{\eeqa}{\end{eqnarray}}
\ncom{\nno}{\nonumber}
\ncom{\non}{\nonumber}
\ncom{\ds}{\displaystyle}
\ncom{\half}{\frac{1}{2}}
\ncom{\mbx}{\makebox{.25cm}}
\ncom{\hs}{\mbox{\hspace{.25cm}}}
\ncom{\rar}{\rightarrow}
\ncom{\Rar}{\Rightarrow}
\ncom{\noin}{\noindent}
\ncom{\bc}{\begin{center}}
\ncom{\ec}{\end{center}}
\ncom{\sz}{\scriptsize}
\ncom{\rf}{\ref}
\ncom{\s}{\sqrt{2}}
\ncom{\sgm}{\sigma}
\ncom{\Sgm}{\Sigma}
\ncom{\psgm}{\sigma^{\prime}}
\ncom{\dt}{\delta}
\ncom{\Dt}{\Delta}
\ncom{\lmd}{\lambda}
\ncom{\Lmd}{\Lambda}
\ncom{\Th}{\Theta}
\ncom{\e}{\eta}
\ncom{\eps}{\epsilon}
\ncom{\pcc}{\stackrel{P}{>}}
\ncom{\lp}{\stackrel{L_{p}}{>}}
\ncom{\dist}{{\rm\,dist}}
\ncom{\sspan}{{\rm\,span}}
\ncom{\re}{{\rm Re\,}}
\ncom{\im}{{\rm Im\,}}
\ncom{\sgn}{{\rm sgn\,}}
\ncom{\ba}{\begin{array}}
\ncom{\ea}{\end{array}}
\ncom{\hone}{\mbox{\hspace{1em}}}
\ncom{\htwo}{\mbox{\hspace{2em}}}
\ncom{\hthree}{\mbox{\hspace{3em}}}
\ncom{\hfour}{\mbox{\hspace{4em}}}
\ncom{\vone}{\vskip 2ex}
\ncom{\vtwo}{\vskip 4ex}
\ncom{\vonee}{\vskip 1.5ex}
\ncom{\vthree}{\vskip 6ex}
\ncom{\vfour}{\vspace*{8ex}}
\ncom{\norm}{\|\;\;\|}
\ncom{\integ}[4]{\int_{#1}^{#2}\,{#3}\,d{#4}}
\ncom{\vspan}[1]{{{\rm\,span}\{ #1 \}}}
\ncom{\dm}[1]{ {\displaystyle{#1} } }
\ncom{\ri}[1]{{#1} \index{#1}}

\newtheorem{remark}{\bf Remark}[section]
\newtheorem{proposition}{Proposition}[section]

\newtheorem{corollary}{Corollary}[section]

\newtheorem{definition}{Definition}[section]
\newtheoremstyle
    {remarkstyle}
    {}
    {11pt}
    {}
    {}
    {\bfseries}
    {:}
    {     }
    {\thmname{#1} \thmnumber{#2} }

\theoremstyle{remarkstyle}

\begin{document}
\newpage

\begin{center}
{\Large \bf Fractional Normal Inverse Gaussian Process}\\
\end{center}

\vone
\begin{center}
{\bf Arun Kumar and P. Vellaisamy}\\
{\it Department of Mathematics\\
Indian Institute of Technology Bombay, Mumbai-400076\\
India}\\
\end{center}
\vtwo
\begin{center}
\noindent{\bf Abstract}
\end{center}
Normal inverse Gaussian (NIG) process was introduced by Barndorff-Nielsen (1997) by subordinating Brownian motion with drift to an inverse Gaussian process.
Increments of NIG process are independent and stationary.
In this paper, we introduce dependence between the increments of NIG process, by subordinating fractional Brownian motion to an 
inverse Gaussian process and call it fractional normal inverse Gaussian (FNIG) process. The basic properties of this process are discussed. 
Its marginal distributions are scale mixtures of normal laws, 
infinitely divisible for the Hurst parameter $1/2\leq H< 1$ and are heavy tailed. First order increments of 
the process are stationary and possess long-range dependence (LRD) property. It is shown that they have persistence of signs LRD property also. 
A generalization of the FNIG process called n-FNIG process is also discussed which allows Hurst parameter $H$ in the interval $(n-1, n)$. Possible 
applications to mathematical finance and hydraulics are also pointed out.
\vone

\section{Introduction}
In recent years, there has been a considerable interest on the study of subordinated processes, since they have wide applications in modeling of financial
time series and hydraulic conductivity in geophysics. Also, these processes have interesting probabilistic properties (Sato (2001)).
A subordinator is a one dimensional 
Levy process with non-decreasing sample paths (a.s.) (Applebaum (2004) p. 49). 
A subordinated process model formulation was originally considered by Mandelbrot and Taylor (1967). 
They considered a one-sided stable law of index $\alpha$ (with $0<\alpha<1)$ for the distribution of the subordinator $\{T(t), t\geq 0\}$ with subordination 
to a Brownian motion $\{B(t), t\geq 0\}$, which gives a 
symmetric stable distribution of index $2\alpha$ for $\{B(T(t)), t\geq 0\}$.
They interpreted the subordinator  $\{T(t),t\geq 0\}$ as trading volume (or number of transactions) upto time t. Clark (1973) used lognormal process 
as a subordinator. Madan and Seneta (1990) considered gamma process as a subordinator and they called the 
subordinated process $\{B(T(t)), t\geq 0\}$ variance gamma (VG) process. Barndorff-Nielsen (1997) introduced normal inverse Gaussian (NIG) process which is
obtained by subordinating Brownian motion with drift to an inverse Gaussian process for modeling of stochastic volatility. Heyde (2005) discussed 
Student Levy process which is obtained by taking inverse gamma process as a subordinator.\\
The marginals of these subordinated processes are generally heavy tailed and are more peaked around the mode. The subordinated processes are useful for
modeling of the returns from financial assets (stocks). Indeed, log-returns from financial time series often exhibits a significant dependence
structure ( Shephard (1995)). 
Subordinated processes, with subordination to Brownian motion have independent increments. 
So, the above processes are limited in nature and may not be suitable for financial data which exhibits a dependence structure.\\
Recently, fractional Brownian motion (FBM) is used instead of Brownian motion for subordination, as its leads to LRD property, 
that is, the covariance function $\gamma_k$ = cov$(X(t), X({t+k}))$ of the stationary time series $\{X(t), t\geq 0\}$ decay so slowly such 
that the series sum $\sum_{k=1}^{\infty} \gamma(k)$ diverges (Beran (1994)). FBM is a natural
generalization of Brownian motion (Mandelbrot (1968)). The FBM $\{B_H(t), t\geq 0\}$ with Hurst parameter $H$ is a centered Gaussian process with 
$B_H(0) = 0$ and with covariance function
\beq\label{covfbm}
E(B_H(t)B_H(s)) = \frac{\sigma^2}{2}[t^{2H} + s^{2H} - |t-s|^{2H}],~~ t,s\geq 0,
\eeq
where $\sigma^2 = $ Var($B_H(1)$). Also, $B_H(t) \sim N(0, \sigma^2t^{2H})$. The FBM has also the following integral representation (Mandelbrot (1968))
\beq \label{fbm}
B_H(t) - B_H(0) = \frac{1}{\Gamma(H+\frac{1}{2})} \Big\{\int_{-\infty}^{0}[(t-u)^{H-1/2}-(-u)^{H-1/2}]dB(u)+
\int_{0}^{t}(t-u)^{H-1/2} dB(u)\Big\}.
\eeq
Its first order increments $\{X(k) = B_H(k+1) - B_H(k)\}$, called fractional Gaussian noise are stationary and exhibits LRD property, when 
$1/2<H<1$.  
Mandelbrot (1997) considered $\{B_H(\theta(t)), t\geq 0\}$, where $\{B_H(t), t\geq 0\}$ is an FBM and the
 activity time process $\{\theta(t), t\geq 0\}$ is a multifractal process with non-decreasing paths and stationary increments. Recently, Kozubowski 
{\it et.al.} (2006) introduced and studied fractional Laplace motion (FLM) by subordinating FBM to a gamma process. They have pointed out the applications
 of FLM to modeling of hydraulic conductivity in geophysics and financial time series. In the same spirit, we subordinate the FBM to in inverse Gaussian
 process to obtain FNIG process and study their probabilistic properties.
\section{FNIG Process}
\setcounter{equation}{0}
The density function of an inverse Gaussian distribution with parameters $a$ and $b$ denoted by $IG(a, b)$, is given by 
\beq
f(x;a,b) = \displaystyle{(2\pi)}^{-1/2} a x^{-3/2}e^{ab-\frac{1}{2}(a^2 x^{-1} + b^2 x)}.
\eeq
The inverse Gaussian (IG) subordinator with parameters $\alpha$ and $\beta$ is defined 
by (Applebaum (2004), p.51)
\beq
G(t) = inf\{s > 0; B(s) + \beta s = \alpha t \},
\eeq 
where $\alpha > 0$ and $\beta \in \mathbb{R}$. Here, $G(t)$ denotes the first time when Brownian motion with drift $\beta$ hits the barrier $\alpha t$. 
Another interpretation of IG subordinator is that it is a Levy process such 
that the increment $G({t+s}) - G(s) $ has inverse Gaussian $IG(\alpha t, \beta)$ distribution. We call $\{G(t), t\geq 0\}$  as IG process with parameters
$\alpha$ and $\beta$. Barndorff-Nielsen (1997) introduced and studied the subordination of Brownian motion with drift to an IG subordinator. Let $Z(t)\sim$ 
IG$(\delta t, \sqrt{\alpha^2 - \beta^2})$, where $0\leq |\beta|\leq\alpha$ and $\delta>0$. Specifically, he studied the process 
$\{N(t), t\geq 0\}$, where $N(t) = B(Z(t)) + \beta Z(t) + \mu t$ and $\mu\in \mathbb{R},$ called NIG process with parameters 
$\alpha, \beta, \mu$ and $\delta$, denoted by NIG$(\alpha, \beta, \mu, \delta)$. 
The process $\{N(t), t\geq 0\}$ is a Levy process and its marginals have heavier tails than Gaussian distribution. Since NIG process is a Levy process 
its increments do not possess LRD property. Our aim here is to introduce dependence between the increments by taking FBM instead of Brownian motion.
Let $\{B_H(t), t\geq 0\}$ be an FBM and 
let $\{G(t), t\geq 0\}$ be an IG process with parameters $\alpha>0$ and $\beta\geq 0$. Also, we assume $H = 1/2$, when $\beta =0.$
We now define a new process $\{X(t), t\geq 0\}$, by 
subordinating FBM to an IG process:
\beq \label{fnig}
\{X(t),t\geq 0\}\stackrel{d} = \{B_H(G(t)), t\geq 0\},
\eeq 
where $\stackrel{d} =$ stands for equality of finite dimensional distributions. 
We call $\{X(t), t\geq 0\}$ fractional normal inverse Gaussian (FNIG) process with parameters $\theta =(\alpha, \beta, \sigma^2, H)$. 
Note that $X(t)$ is a variance mixture of normal distribution with mean zero.
Consequently, for all $t\geq 0$, $X(t)$ can be represented as
\beq\label{a1}
X(t)\stackrel{d} = \sigma G(t)^H Z,
\eeq 
where $Z$ is a standard normal variable and is independent of $G(t)$ having an inverse Gaussian IG$(\alpha t, \beta)$ distribution.
\begin{remark} In the representation (\ref{fnig}), the distribution of $G(t)$ and $B_H(G(t))$ determine each other, i.e.,
\beq
B_H(G(t_1))\stackrel{d} = B_H(G(t_2))\Longleftrightarrow G(t_1) \stackrel{d} = G(t_2).
\eeq
\end{remark}
\noindent The proof is as follows.
\bea
&&B_H(G(t_1))\stackrel{d} = B_H(G(t_2)) ~~\mbox{for some $t_1$ and $t_2$}\nonumber\\
&&\displaystyle\Longleftrightarrow E[e^{iuB_H(G(t_1))}]= E[e^{iuB_H(G(t_2))}]\nonumber\\
&&\displaystyle \Longleftrightarrow  E e^{-\frac{1}{2}u^2\sigma^2G(t_1)^{2H}} = E e^{-\frac{1}{2}u^2\sigma^2G(t_2)^{2H}}\nonumber\\
&&\displaystyle\Longleftrightarrow G(t_1)^{2H}\stackrel{d}= G(t_2)^{2H}\nonumber\\
&&\Longleftrightarrow G(t_1)\stackrel{d} = G(t_2).
\eea 
The step above the last one follows by the uniqueness of the Laplace transform.
\begin{remark}
\noindent Since $B_H(t)\sim N(0, \sigma^2t^{2H})$, we have 
\beq\label{sub}
\displaystyle E\Bigl[\exp{\Bigl(iu\frac{B_H(G(t))}{t^H}\Bigr)}\Bigr] = E\Bigl[\exp\Bigl(-\frac{1}{2}\sigma^2u^2(\frac{G(t)}{t})^{2H}\Bigr)\Bigr].
\eeq 
Also, by the law of large numbers for the subordinator $G(t)$ (see Bertoin (1996), p.92),
\bea
\frac{G(t)}{t} \stackrel{a.s.}\longrightarrow EG(1) = \frac{\alpha}{\beta}, ~\mbox{as}~t\rightarrow\infty.	
\eea
Using the result that, $X_n\stackrel{d}\longrightarrow X $ iff $Ef(X_n)\rightarrow Ef(X)$ for every bounded continuous function $f$ and taking 
$f(x) = \exp(-\frac{1}{2}u^2x^{2H})$, we have from (\ref{sub}) 
\bea
\displaystyle \frac{B_H(G(t))}{t^H}\stackrel{d}\longrightarrow N\Bigl(0, \sigma^2(\frac{\alpha}{\beta})^{2H}\Bigr), ~\mbox{ as}~~t\rightarrow\infty.
\eea
\end{remark}
 
\subsection{Basic Properties of FNIG Process} 
We first obtain the density of one-dimensional distributions.
\begin{proposition} Let $\{X(t), t\geq 0\}$ be an FNIG process with parameter vector $(\alpha,\beta,\sigma^2, H)$. 
Then the density function of $X(t)$ is given by
\beq\label{dens}
\displaystyle f_{X_t}(x) = \frac{\alpha t}{2\pi\sigma}exp(\alpha\beta t)\int_{0}^{\infty}\frac{1}{y^{H+3/2}}exp[-\frac{1}{2}(\frac{x^2}{\sigma^2y^{2H}} + 
\frac{\alpha^2t^2}{y}+\beta^2y)]dy, ~~x\in\mathbb{R},
\eeq 
which is clearly symmetric about origin.
\end{proposition}
\noindent{\bf Proof.} Since
\bea
P(X(t)\leq x) \displaystyle& =& \int_{0}^{\infty}P(B_H(y)\leq x)dF_{G(t)}(y),
\eea
we have
\begin{eqnarray}
\displaystyle f_{X(t)}(x) &=& \int_{0}^{\infty}f_{B_H(y)}(x)dF_{G(t)}(y)\nonumber\\
\displaystyle &=& \int_{0}^{\infty}\frac{1}{\sigma y^H \sqrt{2\pi}}\exp({-\frac{x^2}{2\sigma^2 y^{2H}}})\frac{1}{\sqrt{2\pi}} 
\alpha t \exp(\alpha\beta t)y^{-3/2} \exp[-\frac{1}{2}(\frac{\alpha^2 t^2}{y} + \beta^2 y)]dy\nonumber\\
&=& \displaystyle \frac{\alpha t}{2\pi\sigma}exp(\alpha\beta t)\int_{0}^{\infty}\frac{1}{y^{H+3/2}}exp[-\frac{1}{2}(\frac{x^2}
{\sigma^2y^{2H}} + \frac{\alpha^2t^2}{y}+\beta^2y)]dy, ~~x\in\mathbb{R},
\end{eqnarray}
the result follows.
$\hfill{\Box}$
\begin{remark} When $\beta =0$ and $H = 1/2$,
\bea
\displaystyle f_{X_t}(x) = \frac{\alpha\sigma t}{\pi (x^2+\alpha^2\sigma^2t^2)}, x\in R,
\eea
and so $\{X(t), t\geq 0\}$ represents a Cauchy process (see Feller 1971, p. 348).
\end{remark}
\noindent Let $K_{\nu}(\omega)$ is the modified Bessel function of third kind with 
index $\nu$. An integral representation of $K_{\nu}(\omega)$ (J$\o$rgensen (1982), p. 170) is
\beq\label{bessel}
K_{\nu}(\omega) = \displaystyle\frac{1}{2}\int_{0}^{\infty}x^{\nu-1}e^{-\frac{1}{2}\omega(x+x^{-1})}dx, ~\omega>0.
\eeq 
Barndorff-Nielsen (1997) considered the NIG$(\alpha, \beta, \mu, \delta)$ distribution with density
\bea
f(x;\alpha,\beta,\mu,\delta) = \displaystyle a(\alpha,\beta,\mu,\delta)q(\frac{x-\mu}{\delta})^{-1}K_1\{\delta\alpha q(\frac{x-\mu}{\delta})\}exp(\beta x),
\eea
where
\beq
\displaystyle a(\alpha,\beta,\mu,\delta) = \pi^{-1}\alpha exp(\delta\sqrt(\alpha^2-\beta^2)-\beta\mu) \mbox{,}~ q(x) = \sqrt(1+x^2),
\eeq
and $K_1$ is modified Bessel function of third kind and index 1. When $H=1/2$, we have from (\ref{dens}) 	
\beq
\displaystyle f_{X_t}(x) = \frac{\alpha\beta t}{\pi}exp(\alpha\beta t){(x^2+\alpha^2t^2\sigma^2)}^{-1/2}
K_1(\frac{\beta}{\sigma}\sqrt{x^2+\alpha^2t^2\sigma^2})
\eeq
which is same as the distribution of NIG$(\beta/\sigma, 0, 0, \alpha\sigma t)$.\\  
\noindent We have the following result concerning the moments of FNIG process. 
\begin{proposition}\label{mom1}
If $\{X(t), t\geq 0\}$ is the FNIG process with parameters $\alpha$ and $\beta $, then
\beq\label{moment}
E|X(t)|^q = c_q \sqrt{\frac{2}{\pi}}\alpha (\frac{\alpha}{\beta})^{qH-1/2}t^{qH+1/2}\exp(\alpha\beta t)\displaystyle{ K_{1/2-qH}(\alpha\beta t)},	
\eeq
where $\displaystyle c_q = \sqrt{\frac{2^q}{\pi}}\sigma^q\Gamma(\frac{1+q}{2})$
\end{proposition}
\noindent{\bf Proof.} Using (\ref{a1}), 
\bea
E|X(t)|^q  = \sigma^q E|G(t)|^{qH}E|Z|^q, 
\eea
since $G(t)$ and $Z$ are independent. 
The r-th moment of IG($\alpha t, \beta$) distribution (J$\o$rgensen (1982), p. 13) is given by
\begin{eqnarray}\label{mom}
\mu_r' = \frac{K_{r-1/2}(\alpha \beta t)}{K_{-1/2}(\alpha\beta t)}(\frac{\alpha t}{\beta})^r, ~~ r\in \mathbb{R}.
\end{eqnarray}
Using $K_{1/2}(\omega) = K_{-1/2}(\omega) = \displaystyle\sqrt{\frac{\pi}{2}}
\omega^{-1/2}e^{-\omega}$ (J$\o$rgensen (1982), p. 170), we now obtain 
\bea
EG_t^{qH} = \displaystyle\sqrt{\frac{2}{\pi}}\alpha(\frac{\alpha}{\beta})^{qH-1/2}t^{qH+1/2}exp(\alpha\beta t)K_{qH-1/2}(\alpha\beta t).
\eea
The result follows by using the fact that $E|Z|^q$ = $\displaystyle\sqrt{\frac{2^q}{\pi}}\Gamma(\frac{1+q}{2})$ for a standard normal variate $Z$. 
$\hfill{\Box}$
\begin{remark} 
Note that for large $\omega$, $K_{\nu}(\omega) \sim \displaystyle 
\sqrt{\frac{\pi}{2}}e^{-\omega}\omega^{-1/2}$  (J$\o$rgensen (1982) p. 171)
which shows that, $E|X(t)|^q $ $\displaystyle\sim c_q (\frac{\alpha}{\beta})^{qH}t^{qH}$, as 
$t\rightarrow \infty$. This implies the increments $Y_{\eta}(t) = X(t+\eta) - X(t)$ of the FNIG process behave like the increments of FBM,
as the lag $\eta \rightarrow \infty$.
\end{remark} 
\begin{remark} 
Let $\{X(t), t\geq 0\}$ be the FNIG process with parameter $\theta = (\alpha, \beta, \sigma^2, H)$. 
Using (\ref{moment}), the kurtosis of ${X(t)}$ is given by,
\bea
\beta_2(t)& = &\frac{EX(t)^4}{(EX(t)^2)^2}\nonumber\\
&=& 3\sqrt{\frac{\pi}{2}}e^{-\alpha\beta t}(\alpha\beta t)^{-1/2}\frac{K_{1/2-4H}(\alpha\beta t)}{(K_{1/2-2H}(\alpha \beta t))^2}.
\eea
\end{remark}
\noindent A plot of $\beta_2(t)$ against $t$ shows that $\beta_2$ is a decreasing function in t. This implies that the marginals of $X(t)$
(having $\beta_2 >3,$ for all $t\geq 0$) have heavier tails and a larger degree of peakedness as compared to the Gaussian distribution. Since  
$\beta_2 \rightarrow 3,$ as $t \rightarrow\infty$, the process $X(t)$ becomes Gaussian with an increasing lag. Observe from Figure 2 that the 
 density of the standardized (variance unity) random variable $X(1)$, denoted by $f(x)$ have higher degree of peakedness and heavier tails than standard 
normal density denoted by $g(x)$. Also, the density of $X(1)$ tends to a normal density for large values of $\alpha$. \\
\begin{figure}[h]
\centering{\includegraphics[]{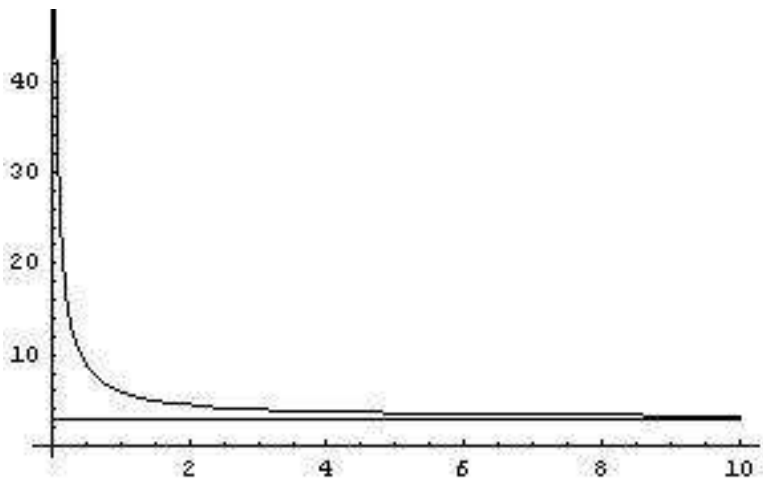}}
\caption{Graph of $\beta_2(t)$ versus t}
\end{figure}

\begin{figure}[h]
\centering{\includegraphics[width= \textwidth]{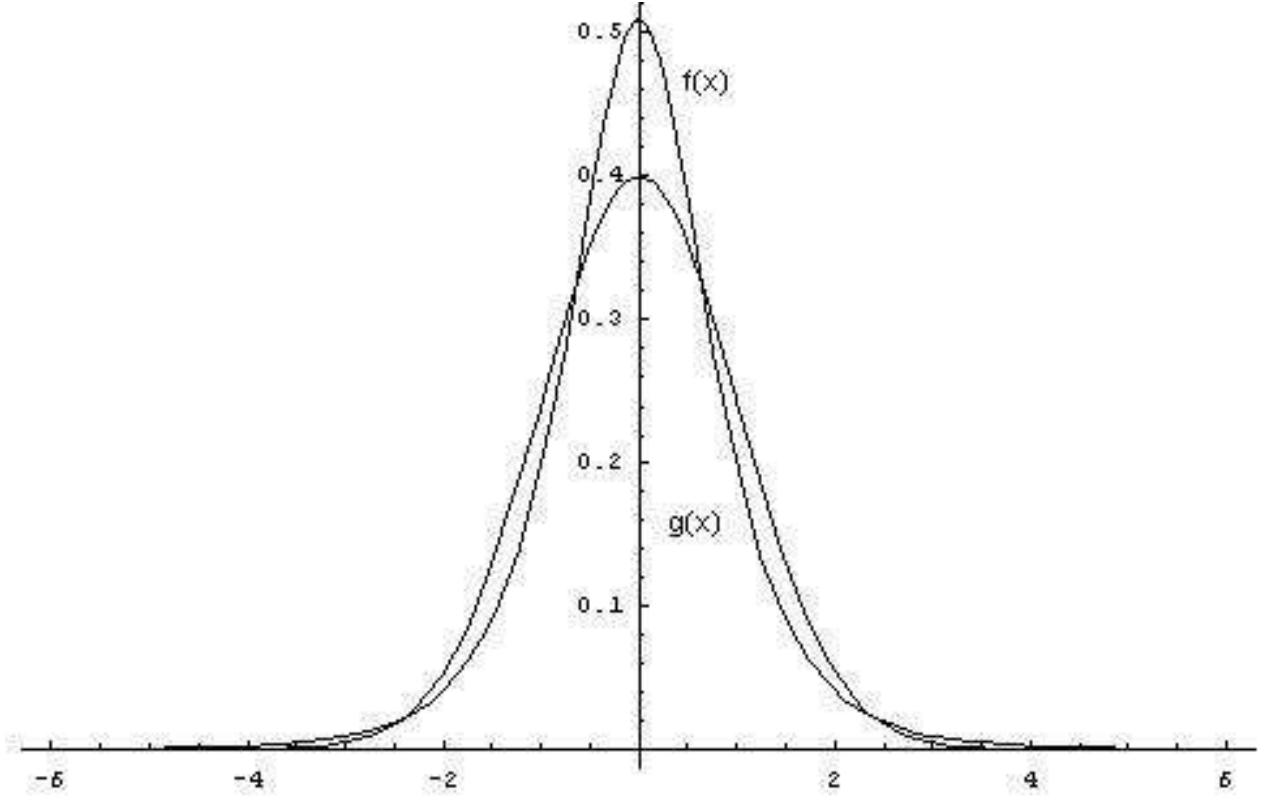}}
\caption{Comparison of pdf of standardized random variable $X(1)$ with parameter $\theta = (3, 1, 1, 0.75)$ and the pdf of $Z\sim N(0,1)$.}
\end{figure}


\subsection{The Covariance Structure}
The following result concerning the covariance function of FNIG process, follows using the fact 
$EX(t)X(s) = \frac{1}{2}[ EX(t)^2 + EX(s)^2 - E(X(t) - X(s))^2 ]$ and using (\ref{moment}). 
\begin{proposition}
The covariance function of FNIG process $\{X(t), t\geq 0\}$ for $t\neq s$, is given by
\begin{eqnarray}
EX(t)X(s) &=& \frac{\sigma^2}{\sqrt{2\pi}}\alpha(\frac{\alpha}{\beta})^{2H-1/2}\Bigl[t^{2H+1/2}\exp(\alpha\beta t)K_{1/2-2H}(\alpha\beta t)\nonumber\\ 
&&\hspace{.6cm} + s^{2H+1/2}\exp(\alpha\beta s)K_{1/2-2H}(\alpha\beta s)\nonumber\\
&&\hspace{.9cm} - |t-s|^{2H+1/2}\exp(\alpha\beta (|t-s|)) K_{1/2-2H}(\alpha\beta(|t-s|))\Bigr].
\label{cov2}
\end{eqnarray}
\end{proposition}
\begin{proposition}\label{aa}
For $\alpha=\beta $,~ $EX(t)X(s)\sim EB_H(t)B_H(s)$ as $t\rightarrow \infty,~ s\rightarrow \infty$ and $|t-s|\rightarrow \infty.$  
\end{proposition}
\noindent{\bf Proof.} Without loss of generality, consider the case 
$0<s<t$, where $s\rightarrow\infty$ with $t-s \rightarrow\infty$. Let $F(s,t) = \displaystyle\frac{EX(t)X(s)}{EB_H(t)B_H(s)}$. We discuss the following
three cases. \\
\noindent {\it Case1.} Let $\displaystyle\frac {s}{t}\rightarrow 1$. Then $\displaystyle(\frac{t-s}{t})\rightarrow 0$, and from (\ref{covfbm}),
 $EB_H(t)B_H(s)\sim \sigma^2t^{2H}$. Hence, using (\ref{cov2}),  
\bea
F(s,t) &\sim& \frac{1}{\sqrt(2\pi)}\alpha[t^{1/2}e^{\alpha^2 t}K_{1/2-2H}(\alpha^2 t) + s^{1/2}e^{\alpha^2 s}
K_{1/2-2H}(\alpha^2 s) \nonumber\\
&&\hspace{.6cm}- (\frac{t-s}{t})^{2H}(t-s)^{1/2}e^{\alpha^2(t-s)K_{1/2-2H}(\alpha^2(t-s))}]\nonumber\\
& \rightarrow & 1. 
\eea
\noindent {\it Case 2.} Suppose that $\displaystyle\frac{s}{t}\rightarrow q$ for some $0<q<1$, so that $\displaystyle(\frac{t-s}{t})\rightarrow (1-q)$, as  
$s\rightarrow\infty.$ In this case, we have
\bea
EB_H(t)B_H(s)\sim \frac{\sigma^2}{2}[t^{2H}(q^{2H}+1-(1-q)^{2H})].
\eea
Hence,
\bea
F(s,t) &\sim& \sqrt\frac{2}{\pi}\alpha\frac{1}{q^{2H}+1-(1-q)^{2H}}[t^{1/2}e^{\alpha^2 t}K_{1/2-2H}(\alpha^2 t)
\nonumber\\
&&\hspace{.6cm}+ q^{2H}s^{1/2}e^{\alpha^2 s}K_{1/2-2H}(\alpha^2 s) - (1-q)^{2H}(t-s)^{1/2}e^{\alpha^2(t-s)}K_{1/2-2H}(\alpha^2(t-s))]\nonumber\\
&\rightarrow& 1.
\eea  
\noindent {\it Case 3.} Suppose that $\displaystyle\frac{s}{t}\rightarrow 0,$ so that $\displaystyle(\frac{t-s}{t}) \rightarrow 1$. In this case,  
\beq\label{eq23}
EB_H(t)B_H(s)\sim\frac{\sigma^2}{2}t^{2H}[1-(1-\frac{s}{t})^{2H}]\sim \sigma^2Hst^{2H-1}. 
\eeq 
An asymptotic expansion for $K_{\nu}(\omega)$ for large $\omega$ (see equation A.9 of J{\o}rgensen (1982), p.171) is
\beq\label{large}
K_{\nu}(\omega) = \sqrt{\frac{\pi}{2}}{\omega}^{-1/2}exp(-\omega)(1 + \frac{\mu-1}{8\omega} + \frac{(\mu-1)(\mu-9)}{2!(8\omega)^2}
+\frac{(\mu-1)(\mu-9)(\mu-25)}{3!(8\omega)^3} +\cdots),
\eeq
where $\mu = 4\nu^2.$  Using the above result, we get $\mu = (1-4H)^2$, and 
\begin{eqnarray}
EX(t)X(s) &=& \frac{\sigma^2}{2}\Bigl[t^{2H}(1 + \frac{\mu-1}{8\alpha^2 t} +
 \frac{(\mu-1)(\mu-9)}{2!(8\alpha^2 t)^2} + \cdots)\nonumber\\ 
&&\hspace{.8cm} +  s^{2H}(1 + \frac{\mu-1}{8\alpha^2 s} + 
\frac{(\mu-1)(\mu-9)}{2!(8\alpha^2 s)^2} + \cdots)\nonumber\\ 
&&\hspace{1.2cm} - (t-s)^{2H}(1 + \frac{\mu-1}{8\alpha^2(t-s)} + \frac{(\mu-1)(\mu-9)}{2!(8\alpha^2(t-s))^2} + \cdots)\Bigr]\nonumber\\
 &=& \frac{\sigma^2}{2}t^{2H}\Bigl[ (1 + \frac{\mu-1}{8\alpha^2 t} + \frac{(\mu-1)(\mu-9)}{2!(8\alpha^2 t)^2} + \cdots) \nonumber\\
&&\hspace{.8cm} +  (\frac{s}{t})^{2H}(1 + \frac{\mu-1}{8\alpha^2 s} + 
\frac{(\mu-1)(\mu-9)}{2!(8\alpha^2 s)^2} + \cdots)\nonumber\\ 
&&\hspace{1.2cm}- (1-\frac{s}{t})^{2H}(1 + \frac{\mu-1}{8\alpha^2(t-s)} + \frac{(\mu-1)(\mu-9)}{2!(8\alpha^2(t-s))^2} + \cdots) \Bigr]\nonumber\\
&\sim&\sigma^2 Hst^{2H-1},
\label{eq34}
\end{eqnarray}
which follows by using binomial series expansion for $\displaystyle(1-\frac{s}{t})^{2H}$, and then allowing $s,t \rightarrow \infty$ such that 
$\displaystyle \frac{s}{t}\rightarrow 0.$ The result now follows from (\ref{eq23}) and (\ref{eq34}).  $\hfill{\Box}$

\subsection{FNIG noise}
Let $\{X(t), t\geq 0\}$ be an FNIG process. Then the first-order increment process 
\beq\label{noise}
\{Y_{\eta}(t), t\geq 0 \}\stackrel{d} = \{ X(t+\eta) - X(t), t\geq 0\}, 
\eeq
for a lag $\eta > 0$ is stationary. 
  
\noindent We have the following result on persistence of signs LRD introduced by Heyde (2002). A process
$\{X(t), t\geq 0\}$ is said to have persistence of signs LRD property if LRD holds for $\{sgn X(t), t\geq 0\}$ (see Heyde (2002)), where
\begin{eqnarray}
sgn(x) = \left\{
\begin{array}{rl}
\displaystyle +1,~~\mbox{if}~~ x>0,\\
0, ~~\mbox{if} ~~x=0,\\
-1, ~~\mbox{if} ~~x<0.
\end{array}\right. 
\end{eqnarray}  
\begin{proposition}
The increment process $\{Y_{\eta}(t), t\geq 0\}$ has persistence of signs LRD property, when $1/2<H<1$.
\end{proposition}
\noindent{\bf Proof.}
Note that $\{sgn Y_{\eta}(t), t\geq 0\}$ has the same finite dimensional distribution as that of $\{sgn(B_H(t+\eta) - B_H(t)), t\geq 0\}$,
which has persistence of signs LRD property (see Proposition 1 of Heyde(2002)). $\hfill{\Box}$
\\\\
We call the process
\beq
\{W_j, j\in\mathbb{N}\} \stackrel{d}= \{Y_{\eta}(\eta(j-1)), j\in\mathbb{N}\}
\eeq
FNIG noise. Next result on the covariance function of FNIG noise directly follows from (\ref{cov2}) and Proposition (\ref{mom1}). 
\begin{proposition}
The covariance function of the process $\{Y_{\eta}(t), t\geq 0\},$ defined by (\ref{noise}), is, for $s\neq t$, given by
\begin{eqnarray}\label{cov1}
EY_{\eta}(t)Y_{\eta}(s)& =&
\frac{\sigma^2}{\sqrt{2\pi}}\alpha(\frac{\alpha}{\beta})^{2H-1/2}\Bigl[|t-s+\eta|^{2H+1/2}\exp(\alpha\beta |t-s+\eta|)K_{1/2-2H}(\alpha\beta |t-s+\eta|) 
\nonumber\\
&&\hspace{.8cm} + |s-t+\eta|^{2H+1/2}\exp(\alpha\beta |s-t+\eta|)K_{1/2-2H}(\alpha\beta |s-t+\eta|)\nonumber\\
&& \hspace{1.2cm}- 2(|t-s|)^{2H+1/2}\exp(\alpha\beta(|t-s|))
K_{1/2-2H}(\alpha\beta(|t-s|))\Bigr],
\end{eqnarray}
and for $s=t,$
\bea
EY_{\eta}(t)Y_{\eta}(s) = E[X(\eta)]^2 = \displaystyle\sigma^2\sqrt{\frac{2}{\pi}}\alpha(\frac{\alpha}{\beta})^{2H-1/2}\eta^{2H+1/2}\exp(\alpha\beta \eta)
K_{1/2-2H}(\alpha\beta \eta). 
\eea
\end{proposition}
\noindent The next result shows the asymptotic behavior of the covariance function of $\{Y_{\eta}(t), t\geq 0\}$.
\begin{proposition} Let $t\geq 0$ be fixed. Then as $s\rightarrow \infty,$
\beq 
EY_{\eta}(t)Y_{\eta}(s)\sim \sigma^2(\frac{\alpha}{\beta})^{2H}\eta^2H(2H-1)(s-t)^{2H-2}.
\eeq
\end{proposition}
\noindent{\bf Proof.} Using the asymptotic expansion for $K_{\nu}(\omega)$ given in (\ref{large}), 
we have,  
\begin{eqnarray}\label{asym1}
EY_{\eta}(t)Y_{\eta}(s) &=& \frac{\sigma^2}{2}(\frac{\alpha}{\beta})^{2H}\Bigl[(s-t-\eta)^{2H}(1 + \frac{(\mu-1)}{8\alpha\beta(s-t-\eta)} +
\frac{(\mu-1)(\mu-9)}{2!(8\alpha\beta(s-t-\eta))^2} +\cdots ) \nonumber\\
&&\hspace{.8cm} + (s-t+\eta)^{2H}(1 + \frac{(\mu-1)}{8\alpha\beta(s-t+\eta)} +
\frac{(\mu-1)(\mu-9)}{2!(8\alpha\beta(s-t+\eta))^2} +\cdots ) \nonumber\\ 
&&\hspace{1.2cm}- 2(s-t)^{2H}(1 + \frac{(\mu-1)}{8\alpha\beta(s-t)} +
\frac{(\mu-1)(\mu-9)}{2!(8\alpha\beta(s-t))^2} +\cdots )\Bigr] \nonumber\\ 
&=& \frac{\sigma^2}{2}(\frac{\alpha}{\beta})^{2H}(s-t)^{2H}\Bigl[(1-\frac{\eta}{s-t})^{2H}(1 + \frac{(\mu-1)}{8\alpha\beta(s-t-\eta)} +
\frac{(\mu-1)(\mu-9)}{2!(8\alpha\beta(s-t-\eta))^2} +\cdots ) \nonumber\\
&&\hspace{.8cm} + (1+\frac{\eta}{s-t})^{2H}(1 + \frac{(\mu-1)}{8\alpha\beta(s-t+\eta)} +
\frac{(\mu-1)(\mu-9)}{2!(8\alpha\beta(s-t+\eta))^2} +\cdots ) \nonumber\\ 
&&\hspace{1.2cm}- 2(1 + \frac{(\mu-1)}{8\alpha\beta(s-t)} +
\frac{(\mu-1)(\mu-9)}{2!(8\alpha\beta(s-t))^2} +\cdots )\Bigr] \nonumber\\
&=& \sigma^2 (\frac{\alpha}{\beta})^{2H}\eta^2 H(2H-1)(s-t)^{2H-2}[1 + O((s-t)^{-1})]. 
\end{eqnarray}
The last equality follows by using binomial series expansion of  $\displaystyle(1+\frac{\eta}{s-t})^{2H}$ and  $\displaystyle(1-\frac{\eta}{s-t})^{2H}$. 
From (\ref{asym1}),
\beq\label{asym}
EY_{\eta}(t)Y_{\eta}(s)\sim\sigma^2(\frac{\alpha}{\beta})^{2H}\eta^2H(2H-1)(s-t)^{2H-2},
\eeq
since $\displaystyle\lim_{s\rightarrow\infty} O((s-t)^{-1}) =0$.
$\hfill{\Box}$
\begin{corollary}
FNIG noise process has LRD property for $1/2 < H<1$. 
\end{corollary}
\noindent{\bf Proof}. The covariance function of $\{W_j, j\geq 1\}$ can easily be obtained from (\ref{cov1}) as
\bea
\gamma(k) = cov(W_j, W_{j+k}) &=& \frac{\sigma^2}{\sqrt{2\pi}}\alpha(\frac{\alpha}{\beta})^{2H-1/2}\Bigl[(\eta(k+1))^{2H+1/2}exp(\alpha\beta\eta(k+1))
K_{1/2-2H}(\alpha\beta\eta(k+1)\\\nonumber
&&\hspace{.8cm}+ (\eta(k-1))^{2H+1/2}exp(\alpha\beta\eta(k-1))K_{1/2-2H} (\alpha\beta\eta(k-1)\\\nonumber
&&\hspace{1.2cm}- 2(\eta k)^{2H+1/2} exp(\alpha\beta\eta k)K_{1/2-2H} (\alpha\beta\eta k)\Bigr], 
\eea
for $k\geq 1.$\\
Further, we have 
\beq\label{lrd}
\gamma(k)\sim \sigma^2(\frac{\alpha}{\beta})^{2H}\eta^{2H}H(2H-1)k^{2H-2},~~\mbox{as}~k\rightarrow\infty.
\eeq
From (\ref{lrd}), we have $\displaystyle \sum_{k=1}^{\infty} \gamma(k) = \infty$, for $1/2<H<1$.  
$\hfill{\Box}$
\begin{remark} Observe that when $\alpha = \beta$ and $\eta = 1$, FNIG noise behaves like fractional Gaussian noise with an 
increasing lag $k$. Also, the covariances 
(and correlations) are positive when $H>1/2$ and negative when $H<1/2$, as in the case of FBM. 
\end{remark} 
\noindent A stochastic process $\{Z(t), t\geq 0 \}$ is said to be
self-similar with Hurst parameter $H\in(0,1)$ if for any $a>0$, we have
\beq
\{Z(at), t\geq 0\} \stackrel{d} = \{a^H Z(t), t\geq 0\}.
\eeq
A stationary process $\{Y(t), t\geq 0\}$ is related to a self-similar process $\{Z(t), t\geq 0\}$ by the Lamperti transformation defined by
\beq\label{lam}
\displaystyle Z(t) = t^H Y(ln(t)), ~~Z(0) = 0, ~~~t>0,
\eeq
(see Lamperti (1962)).
The process $\{Z(t), t\geq 0\}$, is self-similar with exponent $H$. The process $\{Z(t), t\geq 0\}$ may or may not have strictly stationary increments.
If we apply the Lamperti transformation (\ref{lam}) to $Y_{\eta}(t)$ defined in (\ref{noise}), we obtain self-similar process $\{W(t), t\geq 0\}$ 
that do not have strictly stationary increments. In fact, its covariance function for $ s\neq t$, is given by   
\begin{eqnarray}
EW(t)W(s)&=&
\frac{\sigma^2}{\sqrt{2\pi}}\alpha(\frac{\alpha}{\beta})^{2H-1/2}(st)^H\Bigl[|ln(\frac{t}{s})+\eta|^{2H+1/2}\exp(\alpha\beta |ln(\frac{t}{s})+\eta|)
K_{1/2-2H}(\alpha\beta |ln(\frac{t}{s})+\eta|) 
\nonumber\\
&&\hspace{.8cm} + |ln(\frac{t}{s}) - \eta|^{2H+1/2}\exp(\alpha\beta |ln(\frac{t}{s}) - \eta|)K_{1/2-2H}(\alpha\beta |ln(\frac{t}{s}) - \eta|)\nonumber\\
&& \hspace{1.2cm}- 2(|ln(\frac{t}{s})|)^{2H+1/2}\exp(\alpha\beta(|ln(\frac{t}{s})|))
K_{1/2-2H}(\alpha\beta(|ln(\frac{t}{s})|))\Bigr].
\end{eqnarray}
\subsection{Infinite Divisibility}
\noindent The generalized gamma convolutions (g.g.c.) were introduced by Thorin (1978a). They have turned out to be a powerful tool in proving infinite 
divisibility of various distributions (Bondesson (1979)).
A distribution on $\mathbb{R}_+$ is said to be a g.g.c. if it is the weak limit of finite convolutions of gamma distributions 
(see Steutel and Van Harn (2004), p.349). Bondesson (1979) proved that any density on $\mathbb{R_+}$ of the form
\beq\label{ggc}
f(x) =  Cx^{\beta -1}\prod_{j=1}^{M}(1+\sum_{k=1}^{N_j}c_{jk}x^{\alpha_{jk}} )^{-\gamma_j}, ~~~x>0,
\eeq
with all the parameters strictly positive, and $\alpha_{jk}\leq 1$, are g.g.c.'s; consequently all densities (distributions) which 
are weak limit of densities of the form (\ref{ggc}) are g.g.c.'s and thus infinitely divisible. Also, 
if $X$ has a density of the form (\ref{ggc}), the weak limit possibilities included, then the density of $X^q,$ $q\geq 1$, is also a g.g.c. 
(Bondesson (1979)) and hence infinitely divisible. 
\begin{proposition} When $\frac{1}{2}\leq H < 1$, the marginal distributions of FNIG process are infinitely divisible.
\end{proposition}
\noindent {\bf Proof:} It is known that inverse Gaussian distributions or more generally generalized inverse Gaussian (GIG) distributions are g.g.c.'s
(Halgreen (1979)) and hence infinitely divisible.
Using Bondesson (1979) result, we conclude that the stochastic variance $G(t)^{2H}$ has a distribution
which is a g.g.c., when $2H\geq 1 $. Since $\displaystyle (B_H(G(t))|G(t))\sim N(0, \sigma^2 G(t)^{2H})$, and $X(t) = B_H(G(t))$ is a variance mixture of
normal distribution, it is infinitely divisible (Kelker (1971)), when $1/2\leq H <1$.
$\hfill{\Box}$
 

\begin{remark} (i)Kelker (1971) discussed by an example that the variance mixtures of normal distributions can be infinitely divisible, even when the
 stochastic variance is not.
However, we have not been able to establish infinite divisibility of FNIG process for $0<H<1/2$.\\
(ii) Halgreen (1979) showed that variance mixtures of normal distribution are self-decomposable,
 if the stochastic variance is a g.g.c. So, the marginals of FNIG process are self-decomposable, when $1/2\leq H <1$.
\end{remark}
\subsection{Small Deviation Probabilities}
Small deviation probability (in the logarithmic level) for a stochastic process $\{X(t), t\in T\}$ deals with the asymptotic behavior of 
\bea
log P\{\sup_{t\in T}|X(t)|<\epsilon\},~~\mbox{as}~~ \epsilon\rightarrow 0.
\eea
For the inverse Gaussian subordinator $\{G(t), t\geq 0\}$, the Laplace exponent $\Phi$, defined by $\displaystyle E[e^{-xG(t)}] = e^{-t\Phi(x)},$ 
is given by
\bea
\Phi(x) = \alpha[\sqrt{\beta^2+2x} -\beta]. 
\eea
We have 
\bea
\liminf_{x\rightarrow\infty}\frac{\Phi(x)}{ln x}>0.
\eea
The following results on the small deviation probabilities for FNIG process $\{X(t), 0\leq t \leq 1\}$, follows immediately by Theorem 2.1 of 
Linde and Shi (2004).\\
(a) The quenched case (conditional probability $P_{\omega} := P(\cdot|G(t)))$: 
\bea
log P_{\omega} \{||X||_q<\epsilon\} \asymp -\Phi(\epsilon^{-1/H}) ~~\mbox{a.s.}~~\mbox{as}~~\epsilon\rightarrow 0.
\eea
(b) The annealed case (product probability):
\bea
log P \{||X||_q<\epsilon\} \asymp -\Phi(\epsilon^{-1/H}), ~~\epsilon\rightarrow 0,
\eea
where $q\in[1,\infty]$, $f(x)\asymp g(x)$ means, $\displaystyle 0 <\liminf_{x\rightarrow 0}\frac{f(x)}{g(x)}\leq \limsup_{x\rightarrow 0}
\frac{f(x)}{g(x)}<\infty$ and $||X||_q$ is the $L^q$-norm. 

\section {nth-Order Fractional Normal Inverse Gaussian (n-FNIG) Process}
\setcounter{equation}{0}
A generalization of FBM called nth-order FBM, is discussed by Perrin {\it et.al.} (2001).  
This generalization allows the Hurst parameter $H$ in the range $(n-1, n)$. 
Note that n-FBM is obtained by subtracting terms upto the nth-order of the limit expansion of the kernel $(t-u)^{H-1/2}$ from (\ref{fbm}), 
leading to the definition of extended FBM:
\begin{definition}
The nth-order fractional Brownian motion (n-FBM) is defined by
\begin{eqnarray}\label{n-fbm}
B_H^n(t)&=&\frac{1}{\Gamma(H+\frac{1}{2})} \Big\{\int_{-\infty}^{0} \Big[(t-u)^{H-1/2} - (-u)^{H-1/2}- (H-\frac{1}{2})(-u)^{H-3/2}t 
-\cdots\nonumber\\
&&\hspace{.8cm} - (H-\frac{1}{2})\cdots (H-\frac{2n-3}{2})(-u)^{H-n+1/2}
\frac{t^{n-1}}{(n-1)!}\Big]dB(u) \nonumber \\
&&\hspace{1.2cm} +  \int_{0}^{t}(t-u)^{H-1/2} dB(u)\Big \}. 
\end{eqnarray}
\end{definition}
\noindent The covariance function of n-FBM is given by 
\begin{eqnarray}\label{cov} 
G_{H, n}{(t, s)} = (-1)^n\frac{C_H^n}{2}\Big\{|t - s|^{2H} - \sum_{j = 0}^{n - 1} (-1)^j  \left( \begin{array}{c}
2H \\
j \\
\end{array} \right)\Bigl[\Big(\frac{t}{s}\Big)^j|s|^{2H} + \Big(\frac{s}{t}\Big)^j|t|^{2H}\Bigr]\Big\},
\end{eqnarray}
where
\begin{eqnarray}
C_H^n = \frac{1}{\Gamma(2H+1)|sin(\pi H)|}.
\end{eqnarray}
Also, the variance of $B_H^n$ is 
\begin{eqnarray}
V[B_H^n(t)] = C_H^n \left(\begin{array}{c}
2H - 1\\
n - 1\\
\end{array}\right)|t|^{2H} = \sigma_n^2 |t|^{2H} \mbox{(say)}.
\end{eqnarray}
To allow Hurst parameter $H$ in the range $(n-1, n)$, we study nth-order FNIG process by subordinating
inverse Gaussian process with nth-order fractional Brownian motion (n-FBM).
We define n-FNIG process as
\begin{eqnarray}\label{flm} 
\{X^n(t),~t\geq 0\}\stackrel{d}= \{B_H^n(G_t),~t\geq 0\}.
\end{eqnarray}
All the properties of marginals of n-FNIG process follow similar to that of FNIG process by replacing $\sigma$ by $\sigma_n$. 
We can find the covariance function of n-FNIG process with the help of (\ref{cov}) s.t.
\begin{eqnarray}
EB_H^n(G(t))EB_H^n(G(s)) &=& E\Bigr[(-1)^n\frac{C_H^n}{2}\Big\{|G(t) - G(s)|^{2H} - \sum_{j = 0}^{n - 1} (-1)^j  \left( \begin{array}{c}
2H \\
j \\
\end{array} \right)\nonumber\\
&&\hspace{1cm}\times\Bigl[G(t)^j G(s)^{2H-j} + G(s)^j G(t)^{2H-j}\Bigr]\Big\}\Bigr].
\end{eqnarray}
Also, we have 
\begin{eqnarray}\label{eqr}
E G(t)^j G(s)^{2H-j} &=& E\Bigl[\sum_{k=0}^{j}\left(\begin{array}{c}
j\\
k\\
\end{array}\right)(G(t)-G(s))^k G(s)^{2H-k}\Bigr]\nonumber\\
&=&  E\Bigl[\sum_{k=0}^{j}\left(\begin{array}{c}
j\\
k\\
\end{array}\right)(G(t-s))^k G(s)^{2H-k}\Bigr]\nonumber\\
&=& \frac{2}{\pi}\alpha\beta(\frac{\alpha}{\beta})^{2H}\sum_{k=0}^{j}\left(\begin{array}{c}
j\\
k\\
\end{array}\right)\exp{(\alpha\beta t)}K_{k-1/2}(\alpha\beta(|t-s|))\nonumber\\
&&\hspace{.8cm}\times K_{2H-k-1/2}(\alpha \beta s)(|t-s|)^{k+1/2}s^{2H-k+1/2}.
\end{eqnarray}
The last statement follows by using (\ref{mom}). The covariance function of n-FNIG process can be put in explicit form by using (\ref{eqr})
and (\ref{mom}) that is given by
\begin{eqnarray}
EX^n(t)X^n(s) &=& (-1)^n\frac{C_H^n}{2}\Big[\sqrt{\frac{2}{\pi}}\alpha(\frac{\alpha}{\beta})^{2H-1/2}|t-s|^{2H+1/2}\exp(\alpha\beta |t-s|)K_{2H-1/2}
(\alpha\beta |t-s|)\nonumber\\
&&\hspace{.6cm} - \sum_{j = 0}^{n - 1} (-1)^j\left( \begin{array}{c}
2H\\
j \\
\end{array} \right)\Bigl[\frac{2}{\pi}\alpha\beta(\frac{\alpha}{\beta})^{2H}\sum_{k=0}^{j}\left( \begin{array}{c}
j\\
k \\
\end{array} \right)
\Bigl(\exp(\alpha\beta t)K_{k-1/2}(\alpha\beta(|t-s|))\nonumber\\
&&\hspace{1cm} \times K_{2H-k-1/2}(\alpha\beta s)(|t-s|)^{k+1/2}s^{2H-k+1/2}
+ \exp(\alpha\beta s)K_{k-1/2}(\alpha\beta(|t-s|))\nonumber\\
&&\hspace{1.4cm} \times K_{2H-k-1/2}(\alpha\beta t)(|t-s|)^{k+1/2}t^{2H-k+1/2}\Bigr)\Bigr]\Bigr].
\end{eqnarray} 
For $n=1$, the covariance function of n-FNIG process becomes covariance function of FNIG process. Similar to Proposition (\ref{aa}) we can show that
$EX^n(t)X^n(s)\sim EB_H^n(t)B_H^n(s)$ as $s\rightarrow\infty,$ $t\rightarrow\infty$ and $|t-s|\rightarrow\infty.$
\subsection {nth-Order Fractional Normal Inverse Gaussian Noise}
For a function $f$, we denote the kth-order increments by
\begin{eqnarray}
\Delta_l^{(k)}f(x) =\Delta_l^{(k-1)} f(x + l) -\Delta_l^{(k-1)} f(x), ~~k\geq 1.
\end{eqnarray}
Perrin {\it et.al.} (2001) showed that nth-order increments $\{W_H^n(t) = \Delta_l^{(n)}B_H^n(t), t\geq 0\}$ of n-FBM are stationary and 
 hence are called nth-order fractional Gaussian noise.
Covariance function of nth-order fractional Gaussian noise is given by (see Perrin et.al. (2001))
\begin{eqnarray}\label{cnfbm}
cov(W_H^n(t), W_H^n(t+\tau)) = (-1)^n \frac{C_H^n}{2}\sum_{j = -n}^{n}(- 1)^j \left(\begin{array}{c}
2n\\
n + j\\
\end{array}\right)|\tau + jl|^{2H},
\end{eqnarray}
where $\tau >0$.
The nth-order increments of n-FNIG process are stationary and the associated noise is obtained by taking nth-order
increments of the n-FNIG process.
If $\{X^n(t), t\geq 0\}$ is an n-FNIG process, then, for any $\eta > 0$
\begin{eqnarray}
\{Y_{\eta}^n(t), t\geq 0\}\stackrel{d}= \{\Delta_{\eta}^n X^n(t), t\geq 0\}
\end{eqnarray}
is a stationary process. We call the process
\begin{eqnarray}
\{W_j, j\in\mathbb N\}\stackrel{d}= \{Y_{\eta}^n(\eta(j-1)), j\in\mathbb N\}
\end{eqnarray}
an nth-order fractional inverse Gaussian noise (n-FNIGN).
Using (\ref{cnfbm}), we find the covariance function of n-FNIGN process.
\begin{proposition}
Let $\{Y_{\eta}^n(t), t\geq 0\}$ be the n-FNIGN process . Then, for any $t, s \geq 0$, where $s\neq t$, we have 
\bea
E[Y_{\eta}^n(t)Y_{\eta}^n(s)] &=& \frac{C_H^n}{2}{(-1)}^n\displaystyle\sum_{j=-n}^{n}{(-1)}^j  \left( \begin{array}{c}
2n \\
n + j \\
\end{array} \right)\frac{\alpha\tau}{\sqrt(2\pi)}exp(\alpha\beta\tau)\int_{0}^{\infty}{|y+j\eta|}^{2H}y^{-3/2}\nonumber\\	
&&\hspace{1cm}\times exp[-\frac{1}{2}
(\frac{\alpha^2\tau^2}{y} + \beta^2y)]dy.
\eea 
where $\tau = t-s$.
\end{proposition}
\begin{figure}[h]
\centering{\includegraphics[width= \textwidth]{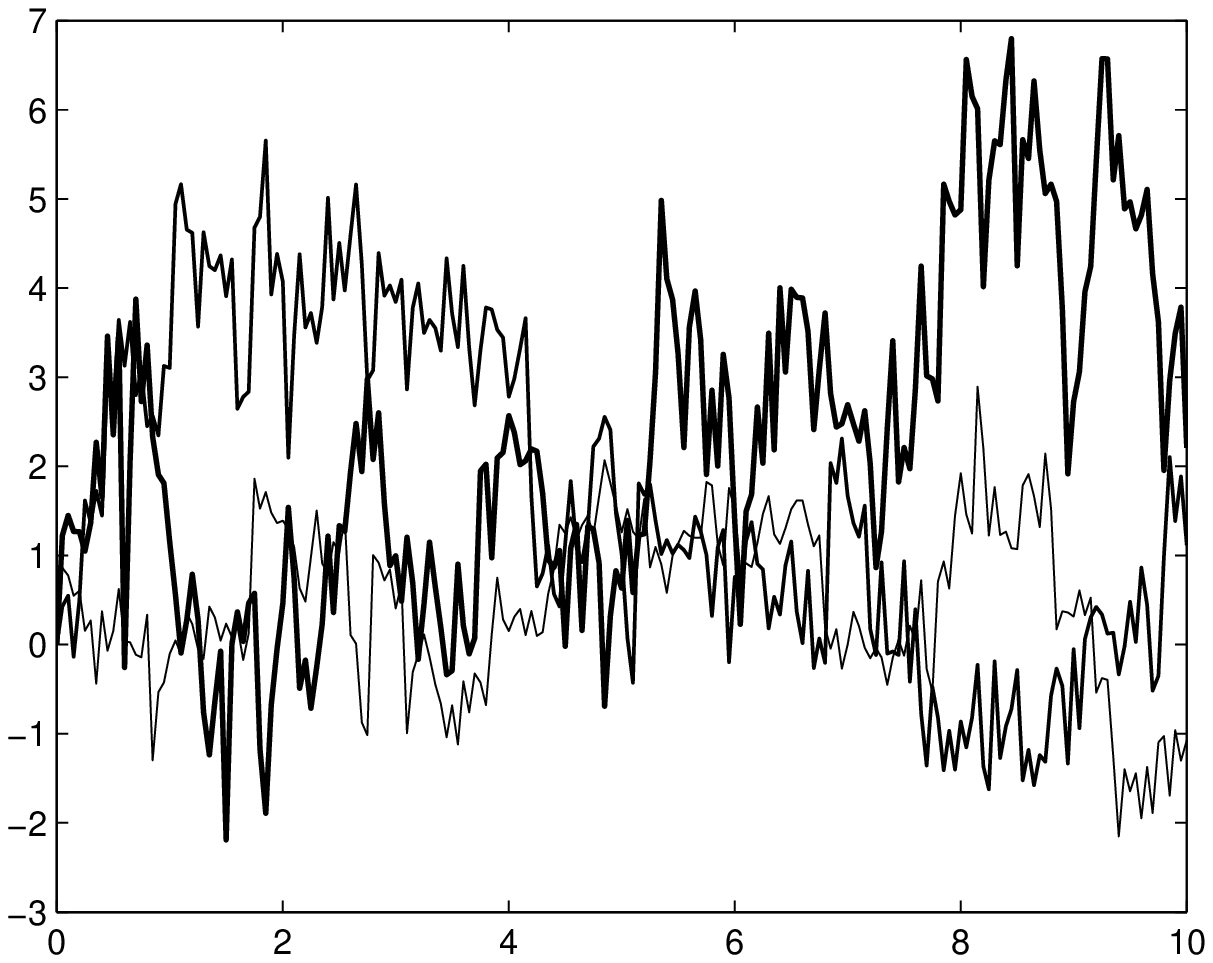}}
\caption{Sample paths of FNIG process for Hurst parameter $H = 0.2$.}
\end{figure}
\begin{figure}[h]
\centering{\includegraphics[width= \textwidth]{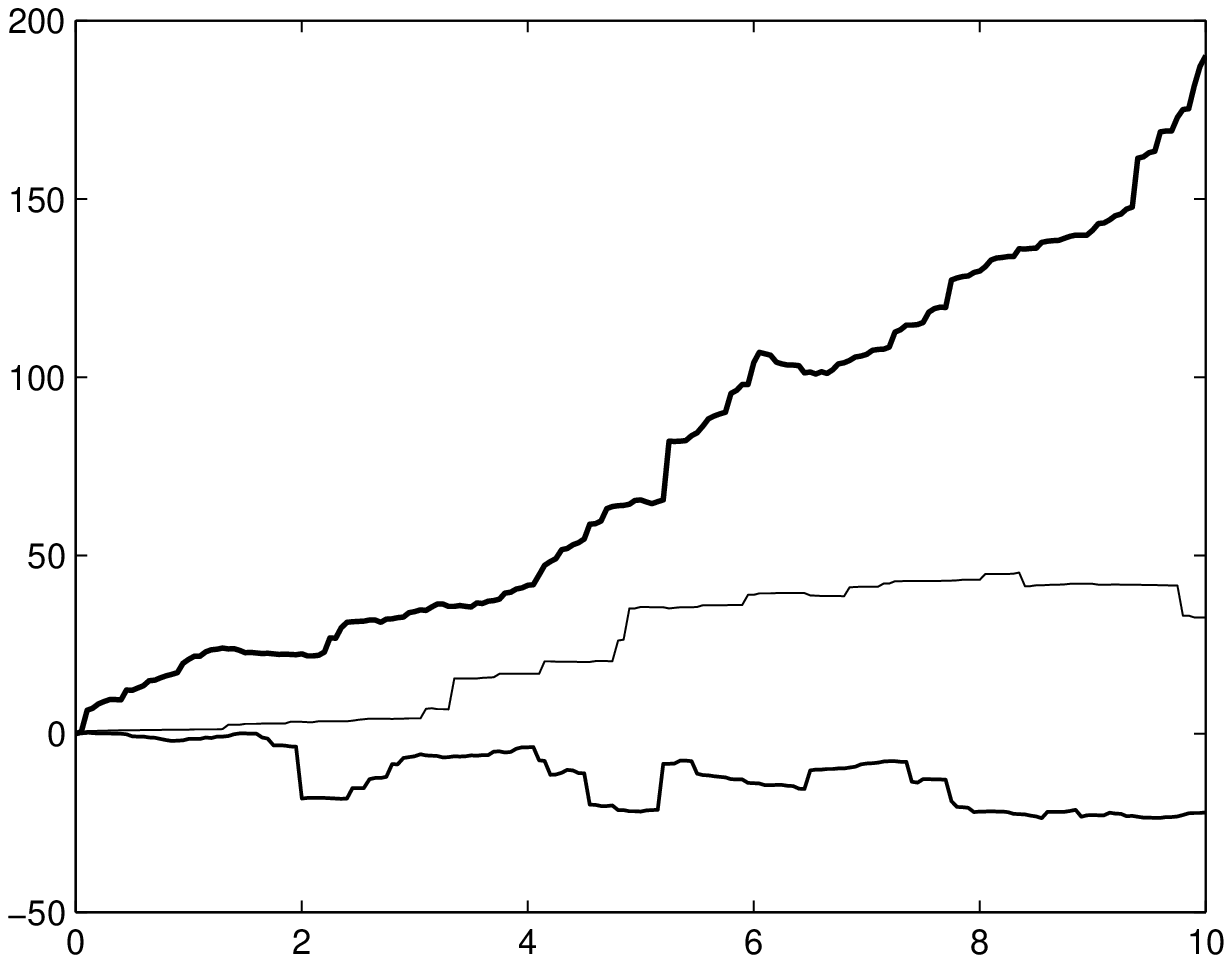}}
\caption{Sample paths of FNIG process for Hurst parameter $H =0.8$.}
\end{figure}

\section{Simulation}
The sample paths of FNIG process can be simulated by subordinating a discretized IG process with fractional Brownian motion. We first simulate a 
discretized IG process on equally spaced intervals. Then we simulate $X_k = B_H(G(t_k))$, which conditionally on the values of $G(t_k)$ is a second order
Gaussian process with the explicit covariance function of FBM evaluated at the values of the IG process. Using Choleski decomposition $LL^T$ of the 
covariance function of the random vector $B_H(G(t_k))|G(t_k), 1\leq k\leq n$. Then $X = (X_1, X_2,\cdots,X_n) = LZ$, where $Z\sim N_n(0, I_n)$,  
a multivariate normal vector. 
The algorithm is given below.
\begin{itemize}
\item [{(i)}] Choose a finite interval $[0, t],$ say our aim to have simulated values at $t_0=0, t_1= t/n,\cdots, t_{n-1}=(n-1)t/n, t_n=t$.
\item[{(ii)}] Generate n IG variates $(G_1, G_2,\cdots,G_n)$, where $G_i = G(t_i)$, by using the algorithm given by Michael, Schucany and Hass 
(see Cont and Tankov (2004), p.182). 
\item[{(iii)}] Form the matrix $\Gamma_{n\times n}=\Bigl(R(G_i,\, G_j)\Bigr)_{n\times n}$, where 
$R(t,s) = \frac{1}{2}(s^{2H} + t^{2H} - |t - s|^ {2H})$.
\item[{(iv)}] Using Cholesky decomposition find $L$ such that $\Gamma_{n\times n}=LL^T$.
\item[{(v)}] Generate a n-dimensional standard normal vector $Z$.
\item[{(iv)}] Compute $X = (X_1, X_2, \cdots, X_n) = LZ$. Then $X_1, X_2, \cdots, X_n$ denotes n-simulated values from FNIG process.
\end{itemize}

\noindent Figures 3 and 4 give the simulated sample paths of FNIG process for the case $\alpha = \beta$, and for the choices 
$\alpha = .1$ (thin lines), $\alpha = 1$ (medium thick lines) and $\alpha =10$ (thick lines). 
Note that sample paths for the case $H = 0.8$ (Figure 4) have less fluctuations, compared to the case when $H=0.2$ (Figure 3). This is because of the 
persistence of signs LRD property of the FNIG process for $1/2<H<1$.   

\section{Applications}
Recently, subordinated processes have been found useful models in finance and geophysics. Most common processes in literature include 
variance gamma (Madan and Seneta (1990)), normal inverse Gaussian (Barndorff-Nielsen (1997)), FATGBM model (Heyde (1999)), 
multifractal model (Mandelbrot {\it et.al.} (1997)), Student processes (Heyde and Leonenko (2005)) and FLM 
(Kozubowski {\it et.al.} (2006)).  
Let $P(t)$ denote the value of a stock at time $t$. Heyde (1999) 
considered the model 
\beq\label{heyde}
\displaystyle P(t) = P(0) \exp({\mu t + \sigma B(T(t))}), 
\eeq
where $\{T(t), t\geq 0\}$ is an increasing process with LRD property and having stationary increments and heavy tails. 
Also, $\{B(t), t\geq 0\}$ is independent of $\{T(t), t\geq 0\}$. 
Then the log-return 
\bea
R(t) = ln\Bigl(\frac{P(t)}{P({t-1})}\Bigr)
\eea
is an uncorrelated process having heavy tails. Mandelbrot et.al. (1997) considered the model
\bea
\displaystyle P(t) = P(0) \exp({B_H(\theta(t))}), 
\eea
where $\{\theta(t), t\geq 0\}$ is a multifractal process with increasing sample paths and stationary increments. It is known that 
NIG distributions are more suitable than Gaussian or generalized hyperbolic distributions (see Barndroff-Nielsen(1997)). 
Shephard (1995) observed that log-returns $R(t)$ often exhibits dependence
structure. Also, the autocorrelations of asset returns are often insignificant, but this is not true for 
small time scales data (Cont and Tankov (2004)). Also, Heyde (2002) observed that 
the log-returns for a  financial time series have persistence of signs LRD property, i.e., the process has the tendency to change signs more rarely than 
under independence or weak dependence. We have shown that FNIG process has LRD property and also have persistence of signs LRD property.
Also, FNIG process is a non-stationary process with stationary 
increments, and reduces to Gaussian process, for sufficiently large lags, and hence could be useful for large-scale models based on fractional Brownian 
motion and non-Gaussian behavior on smaller scales. Processes with these properties were used by Molz and Bowman (1993) and Painter (1996) for analyzing hydraulic
conductivity measurements. Also, Meerschaert {\it et.al.} (2004) used FLM to study the hydraulic conductivity measurements. 
Hence, we believe that FNIG process could serve a useful model in the area of finance, biology, geophysics and etc., 
where the dependence structure of increments and heavy tailedness of 
marginals play an important role. 
\vone
\noindent{\bf Acknowledgment} The first author wishes to thank Council of Scientific and Industrial Research (CSIR), India, for the award of a research 
fellowship. 
\vone
\noindent {\bf \Large References}

\vone
\noindent

\begin{namelist}{xxx}
\item{}  Applebaum, D. (2004). {\it Levy Processes and Stochastic Calculus} (Cambridge Studies in Advanced Mathematics). { Cambridge University Press}, 
Cambridge. 
\item{}  Barndorff-Nielsen, O. E. (1997). Normal inverse Gaussian distributions and stochastic volatility modeling. 
{\it Scand. J. Statist.}, {\bf 24}, 1-13.
\item{}  Beran, J. (1994). {\it Statistics for Long-Memory Processes}. {Chapman \& Hall}, New York.
\item{}  Bertoin, J. (1996). {\it Levy Processes}. {Cambridge University Press}, Cambridge.
\item{}  Bondesson, L. (1979). A general result on infinite divisibility. {\it Ann. Probab.} {\bf 7}(6), 965--979.
\item{}  Clark, P.K. (1973). A subordinated process model with finite variance for speculative prices. {\it Econometrica}, {\bf 41}, 135-155.
\item{}  Cont, R. and Tankov, P. (2004). {\it Financial Modeling with Jump Processes}. { Chapman \& Hall CRC Press}, Boca Raton. 
\item{}  Devroye, L. (1986). {\it Nonuniform Random Variate Generation}. {Springr}, New York. 
\item{}  Feller, W. (1971). {\it Introduction to Probability Theory and its Applications. Vol. $II$.} {John Wiley}, New York.
\item{}  Halgreen, C. (1979). Self-decomposability of the generalized inverse Gaussian and hyperbolic distributions. {\it Z. Wahrsch. Verw.
 Gebiete.} {\bf 47}, 13--17.
\item{}  Heyde, C.C. (1999). A risky asset model with strong dependence through fractal activity time. {\it J. Appl. Probab.}
{\bf 34}(4), 1234-1239.
\item{}  Heyde, C.C. (2002). On modes of long-range dependence. {\it J. Appl. Probab.} {\bf 39}, 882-888.
\item{}  Heyde, C.C. and Leonenko N. N. (2005). Student processes. {\it Adv. Appl. Probab.} {\bf 37}, 342-365.
\item{}  Kelker, D. (1971). Infinite divisibility and variance mixtures of the normal distribution. {\it Ann. Math. Statist.} {\bf 42}, 802-808.
\item{}  Kozubowski, T.J., Meerschaert, M.M., and Podgorski, K. (2006). Fractional Laplace motion. 
{\it Adv. Appl. Prob.} {\bf 38}, 451-464. 
\item{}  Lamperti, J. W. (1962). Semi-stable stochastic processes . {\it Trans. Amer. Math. Soc.} {\bf 104}, 62-78.
\item{}  Linde, W. and Shi, Z. (2004). Evaluating the small deviation probabilities for subordinated Levy processes.
 {\it Stochastic Process. Appl.} {\bf 113}, 273-287.
\item{} Madan, D.B., Carr, P. and Chang, E.C. (1998). The variance gamma process and option pricing. {\it Europ. Finance Rev.}
{\bf 2}, 74-105. 
\item{} Madan, D.B. and Seneta, E. (1990) The variance gamma (V.G) model for share markets returns. {\it J. Business}
{\bf 63}, 511-524.
\item{} Mandelbrot, B.B. (2001). Scaling in financial prices: I. Tails and dependence. {\it Quant. Finance} {\bf 1}, 113-123.
\item{} Mandelbrot, B.B., Fisher, A. and Calvet, L. (1997). A multifractal model of asset returns. {Cowles Foundation Discussion Paper No. 1164}.
\item{} Mandelbrot, B.B. and Taylor, H. (1967). On the distribution of stock price differences. {\it Operat. Res.} {\bf 15}, 1057-1062. 
\item{} Mandelbrot, B.B and Van Ness, J.W. (1968). Fractional Brownian motion, fractional noises and applications.
 {\it SIAM Rev.} {\bf 10}, 422-438.
\item {} Meerschaert, M.M., Kozubowski, T.J., Molz, F.J., and Lu, S. (2004). Fractional
Laplace model for hydraulic conductivity. {\it  Geophys. Res. Lett.} {\bf 31}, p. L08501.
\item{} Molz, F. J., and Bowman, G.K. (1993). A fractal-based stochastic interpolation scheme in subsurface hydrology. {\it Water Resour. Res.} 
{\bf 32}, 1183-1195.
\item{} Painter, S. (1996). Evidence for non-Gaussian scaling behavior in heterogeneous sedimentary formations. {\it Water Resour. Res.} {\bf 32}, 1183-1195.
\item{} Perrin, E., Harba, R., Berzin-Joseph, C., Iribarren, I. and Bonami, A. (2001). nth-Order fractional Brownian 
motion and fractional Gaussian noises. {\it IEEE Trans. Sig. Proc.} {\bf 49}, 1049-1059. 
\item{}  Sato, K. (2001). Subordination and self-decomposability. {\it Statist. Probab. Lett.}  {\bf 54} (3), 317--324.
\item{}  Shephard, N. (1995). Statistical aspects of ARCH and stochastic volatility. {\it In Time series models}. {\it In econometrics, finance and others 
fields} (eds D.R. Cox, D.V. Hinkley \& O.E. Barndorff-Nielsen), 1-67. Chapman \& Hall, London.
\item{}  Steutel, F.W. and Van Harn, K. (2004). {\it Infinite Divisibility of Probability Distributions on the Real Line.}
 Marcel Dekker, New York.
\item{}  Thorin, O. (1978a). An extension of the notion of a generalized gamma convolution. {\it Scand. Actur. J.}, 141-149.
\end{namelist}{}
\end{document}